\documentclass[final]{siamltex}

\usepackage[intlimits]{amsmath}
\usepackage{amssymb}
\usepackage{graphicx}
\usepackage{cancel}
\usepackage{subfigure}

\title{Spontaneous oscillations  in simple fluid networks}
\author{Nathaniel J. Karst\footnotemark[1]\ \and Brian D. Storey\footnotemark[2]\ \and John B. Geddes\footnotemark[2]\ }

\begin{document}
\maketitle

\renewcommand{\thefootnote}{\fnsymbol{footnote}}
\footnotetext[1]{Babson College, Babson Park MA 02457 }
\footnotetext[2]{Olin College, Needham MA 02492 }
\renewcommand{\thefootnote}{\arabic{footnote}}

\begin{abstract}
Nonlinear phenomena including multiple equilibria and spontaneous oscillations are common in fluid networks containing either multiple phases or constituent flows. In many systems, such behavior might be attributed to the complicated geometry of the network, the complex rheology of the constituent fluids, or, in the case of microvascular blood flow, biological control. In this paper we investigate two examples of a simple three-node fluid network containing two miscible Newtonian fluids of differing viscosities, the first modeling microvascular blood flow and the second modeling stratified laminar flow.  We use a combination of analytic and numerical techniques to identify and track saddle-node and Hopf bifurcations through the large parameter space. In both models, we document sustained spontaneous oscillations and, for an experimentally relevant example of parameter analysis, investigate the sensitivity of these oscillations to changes in the viscosity contrast between the constituent fluids and the inlet flow rates. For the case of stratified laminar flow, we detail a physically realizable set of network parameters that exhibit rich dynamics.  The tools and results developed here are general and could be applied to other physical systems.
\end{abstract}

\section{Introduction}

A classic  problem in the field of hydraulics is  determining the distribution of flow rates and pressures inside a given piping network for fixed inlet conditions. Many practical fluid networks such as municipal water delivery have turbulent flow and thus a nonlinear resistance making their analytical solution difficult. In 1936, a structural engineer named Hardy Cross revolutionized the  analysis of hydraulic networks by developing a systematic iterative method by which one could reliably solve nonlinear network problems by hand calculation \cite{Cross1936}.

While analysis of such hydraulic networks is now considered routine with computer techniques, the problem can once again become intractable if one considers networks filled with a fluid comprised of multiple phases or constituents. Analysis of such networks is of interest because in a number of application it has been observed that the phase distribution within the network may exhibit unsteady or non-unique flow. Such heterogeneous distribution of phase within network flows has been studied  at a variety of scales. At the micro-scale, the flow of droplets or bubbles through microfluidic networks can demonstrate bistabilty and spontaneous oscillations~\cite{Jousse2006,Schindler2008,Fuerstman2007,Prakash2007, Joanicot2005}. These nonlinearities  have been exploited by researchers who have demonstrated microfluidic memory, logic, and control devices~\cite{Fuerstman2007,Prakash2007, Joanicot2005}. On the macro-scale,  models of magma flow with either temperature-dependent viscosity~\cite{Helfrich1995}  or volatile-dependent viscosity~\cite{Wylie1999} have shown the existence of multiple solutions on the pressure-flow curve which can lead to spontaneous oscillations.

Another network that can exhibit complex behavior is microvascular blood flow. Nobel prize winner August Krogh noted the heterogeneity of blood flow in the webbed feet of frogs in the early 1920's~\cite{Krogh:1921aa}. In the {\it Anatomy and Physiology of Capillaries} he wrote~\cite{Krogh:1922aa}
\begin{quote}
In single capillaries the flow may become retarded or accelerated from no visible cause; in capillary anastomoses the direction of flow may change from time to time.
\end{quote}
Numerous researchers have confirmed these observations over the years. The heterogeneous distribution of red blood cells in microvascular blood flow is often interpreted as evidence of biological control. If the flow in a branch increases, it is assumed that the diameter of the branch responds in order to auto-regulate the flow. Vasomotion has often been assumed to be the cause for oscillations in the micro-circulation~\cite{rodgers}. While the importance of vasomotion cannot be denied, there is significant evidence that fluctuations in cell distributions in microvascular networks can be due to inherent instabilities~\cite{Kiani:1994aa,Carr:2000aa}.

There are two fundamental phenomena in two-phase  flow networks which differ from their single phase counterparts and lead to complicated behavior. The first effect is that the effective viscosity or flow resistance in a single pipe is often a nonlinear function of the fraction of the different fluids in the pipe. The second effect is that in two fluid systems, it is commonly observed that the phase fraction after a diverging node  is different in the two downstream branches. In 2007 we (JBG and NJK) proved  that if the viscosity is a nonlinear function of fluid fraction then multiple stable equilibrium states may exist~\cite{Geddes:2007}. Further we proved that  both phase separation at a node and nonlinear viscosity can lead to the emergence of spontaneous oscillations. In recent experiments,  we (JBG and BDS) have demonstrated some of these predictions experimentally in simple networks involving two Newtonian fluids of different viscosity~\cite{geddes2010a,karst2013}. In one set of experiments we demonstrated bistability via nonlinear resistance~\cite{geddes2010a} and in the other bistability via phase separation~\cite{karst2013}. These experiments showed that multiple equilibria in networks is possible without fluids with complex rheology.

While our experiments involve simple fluids in a controlled laboratory setting, it is expected that these results may be generalized and found in numerous natural and man-made systems. Phase separation at a single node exists in numerous  fluid systems and has been widely studied in different contexts. In microvascular blood flow, Krogh introduced the term ``plasma skimming''  in order to explain the disproportionate distribution of red blood cells observed at single branch bifurcations {\it in vivo}~\cite{Krogh:1921aa}.  Numerous authors have demonstrated plasma skimming {\it in vitro} and {\it in vivo} and developed simple empirical models to describe the effect~\cite{Bugliarello:1964aa,Chien:1985aa,Dellimore:1983aa,Fenton:1985aa,Klitzman:1982aa,Pries:1989aa}. Another widely studied example of phase distribution at a single node is gas-liquid two-phase flow which has important technological applications in power and process industries. In many process applications phase maldistribution can have detrimental consequences for downstream equipment~\cite{lahey1986}, while in some cases the phenomenon is exploited to build simple phase separators~\cite{Azzopardi1993}. Extensive experimental work on gas-liquid flow has been conducted over the past 50 years~\cite{Azzopardi1999,Azzopardi1994,lahey1986}. In applications for the  process and petroleum industry, phase separation in liquid-liquid flows are  less well-studied though several recent papers  have emerged~\cite{yang2006,yang2007,wang2008}. The impact of phase maldistribution in two-phase flow has been  shown to impact network flows in  refrigeration systems~\cite{hongliang2009} and solar power systems~\cite{minzer2006}.

While the behavior at a single node has been well-studied experimentally in the applications noted above, systematic analysis of networks with two-phase flow have received less attention. The most widely studied network is the microvascular one for which the first modeling effort for dynamics  was  conducted by Kiani {\it et al.}~\cite{Kiani:1994aa}.  In 1994 they  conducted a direct simulation of 400 vessels and found oscillations in the flow. In 2000 Carr and Lecoin  found oscillations in  networks with  fifteen vessels~\cite{Carr:2000aa}. They found evidence of Hopf bifurcations and limit cycles, but were unable to determine which parameters controlled the dynamics. In an attempt to understand the parameters that lead to spontaneous oscillations in microvascular flows, Geddes \textit{et al.} performed a complete analysis of the flow-driven 2-node network (one inlet, a loop, and one outlet) in 2007~\cite{Geddes:2007}. While this network can exhibit oscillations in theory, they do not exist for realistic physical parameters. Several other groups have since studied the problem of oscillations in microvascular networks, and a coherent picture is beginning to emerge~\cite{Pop:2007aa,Shevkoplyas, Owen, Pozrikidis}.

In the context of microvascular flow we now know that networks with 2 vessels can exhibit spontaneous oscillations  for unrealistic physical parameters while networks with 15 vessels can oscillate for realistic parameters~\cite{Carr:2000aa,Geddes:2007}. It is unknown at what level of network complexity  oscillations can emerge and what parameters govern their existence. While the 2-node network has been fully characterized theoretically, full descriptions of more complicated networks becomes difficult. While we have studied the equilibrium properties of the 3-node network (two inlets, a loop, and one outlet) theoretically and experimentally in prior work~\cite{geddes2010a,karst2013}, we had no systematic method to understand the stability  other than through direct simulation. While we did not predict the existence of oscillations for the parameters relevant to our experiments, with no systematic method to analyze the stability and the large number of parameters it is impossible to rule out the emergence of spontaneous oscillations.

In this paper we develop a methodology for finding and tracking  Hopf bifurcations through continuation. This development is critical due to the large parameter space of the problem. We find that our analytical methods are in perfect agreement with direct numerical simulations, validating the methodology. Using our methods we develop phase diagrams that show a rich set of dynamics including multiple frequency oscillations and co-existing limit cycles.   The details of these predictions depend sensitively on the constitutive laws for the fluids in the network and the phase separation at a single diverging node. However, our methodology is general and may be applied to any two-phase flow network system.

\section{Three node model}

\begin{figure}[h!]
\begin{center}
\includegraphics[width=55mm]{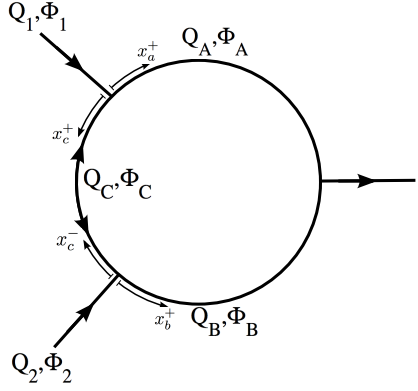}
\end{center}
\caption{Schematic of a configuration of the three node network. Inlet 1 and inlet 2 supply fluids of volume fraction $\Phi_1$ and $\Phi_2$ at controlled flow rates $Q_1$ and $Q_2$. While the flow in vessels $A$ and $B$ is always from left to right in this figure, the flow in vessel $C$ can be up or down depending on the state of the network, with $Q_C > 0$ representing downward flow and $Q_C < 0$ representing upward flow. }\label{fig:schematic}
\end{figure}

The physical setup is shown in Figure \ref{fig:schematic}. The network has two flow controlled inlets, each of which contains a fluid comprised of two separate phases, $\alpha$ and $\beta$. The two phases are two fluids which have different viscosities and remain distinct at least up to the inlet of the network. Without loss of generality, we  assume that $\beta$ is the more viscous fluid. Locally at a point along the tube we  define the local volume fraction as $\Phi = Q_\beta/(Q_\alpha + Q_\beta)$, where $Q$ is the volumetric flow rate of each phase. In an experiment the volume fraction in the two inlets would be set upstream by flow controlled pumps attached to reservoirs of fluids $\alpha$ and $\beta$. In the case of blood flow the $\alpha$ phase is plasma, the $\beta$ phase is red blood cells and the volume fraction is the hematocrit. While blood is not really comprised of two continuous  fluid phases, such a model is commonly used in numerical simulations or laboratory experiments \cite{Popel:2005aa}.

The basic network model based on fundamental conservation principles will be developed in the next section. However, to close the model we require two constitutive laws which depend on the details of the fluid system and the network geometry; i)  the effective viscosity  as a function of volume fraction and ii) the phase separation rule for a single node. The details encoded in these two constitutive laws play a dramatic role in the eventual behavior of the network~\cite{karst2013}.

We assume laminar flow in cylindrical tubes where the hydraulic resistance is proportional to the viscosity  of the fluid mixture. Since we have have a two-phase flow we can compute an effective viscosity, $\mu$, which is a function of not only the two fluids involved but their geometrical arrangement in the tube. The effective viscosity can be expressed in terms of the viscosity of the less-viscous phase, $\mu_\alpha$, and a relative viscosity; $\mu = \mu_\alpha \mu_{\mathrm{rel}}$. Simple Newtonian  fluids  approximately follow a nonlinear Arrhenhius law when they are well mixed,
\begin{equation}
\mu_{\mathrm{rel}} = \left(\frac{\mu_\beta}{\mu_\alpha} \right)^\Phi.
\label{eq:viscosity}
\end{equation}
Here $\mu_\alpha$ and $\mu_\beta$ are the viscosities of the individual phases, and $\mu_\beta /\mu_\alpha$ is the viscosity contrast.

Different viscosity laws exist for different physical manifestations rather than complete mixing. For Newtonian fluids that remain stratified in a circular tube as separate phases, the relative viscosity follows a relationship which can be readily computed though no simple analytical form exists~\cite{geddes2010a,Gemmell1962}. Another common physical configuration  is a core annular flow where the viscous fluid assumes a cylindrical core which is lubricated by an annulus of less viscous fluid in a cylindrical tube~ \cite{Joseph1997}. For the example of microvascular blood flow the rheology is more complicated, however Pries \emph{et al.}~\cite{Pries:1992aa} compiled a database of viscosity measurements in tubes with a range of diameters and  hematocrits.  While the exact form of the above viscosity laws all differ, the important fact is that they are all nonlinear functions of the volume fraction which is a key feature for networks to exhibit  multiple equilibrium states and spontaneous oscillations~\cite{Geddes:2007}. Throughout this work we will assume for convenience that the effective viscosity is determined by Equation \ref{eq:viscosity}.

\begin{figure}[!t]
    \begin{subfigure}[]
                \centering
                \includegraphics[width=65mm]{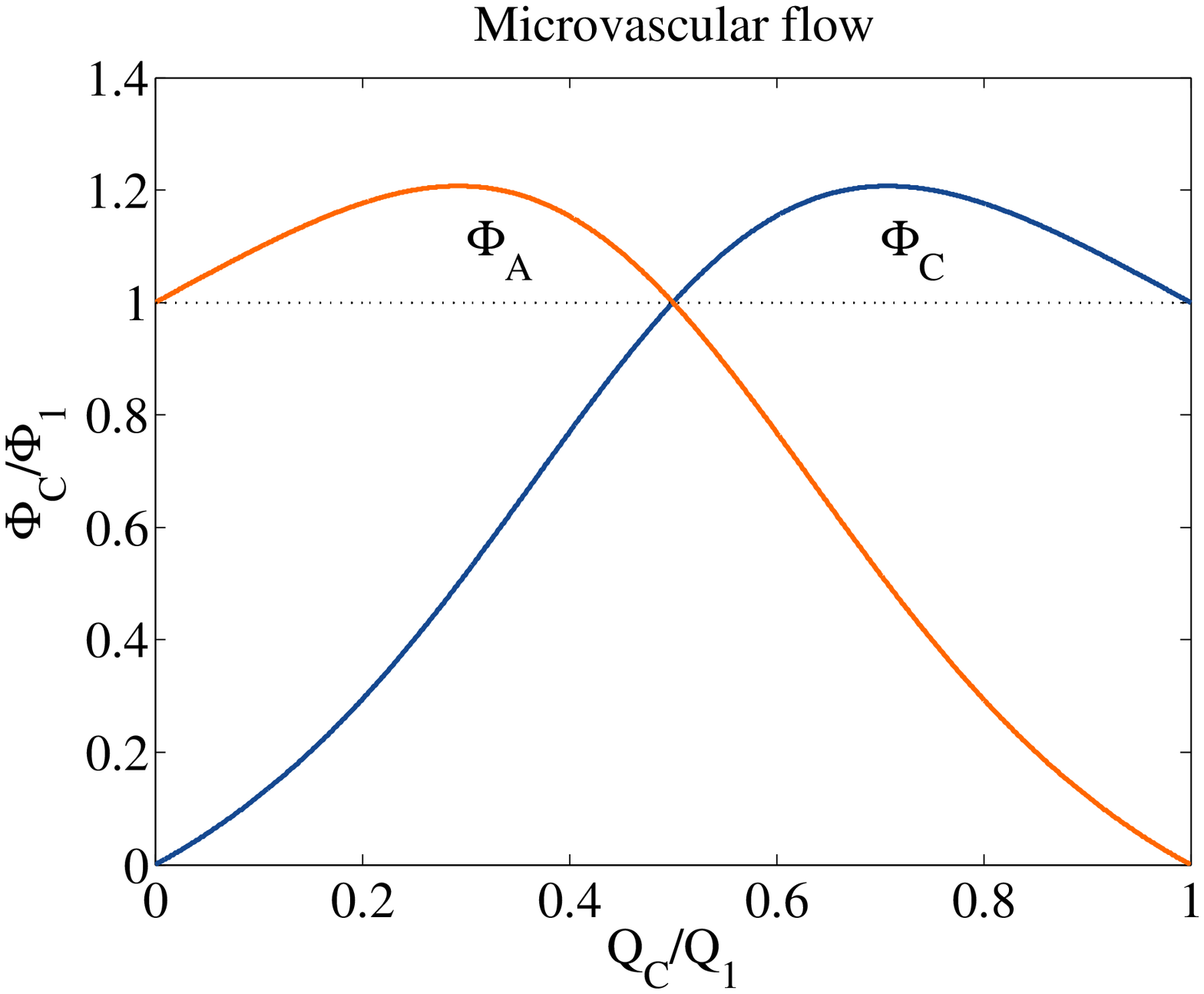}
        \end{subfigure}
        \begin{subfigure}[]
                \centering
                \includegraphics[width=65mm]{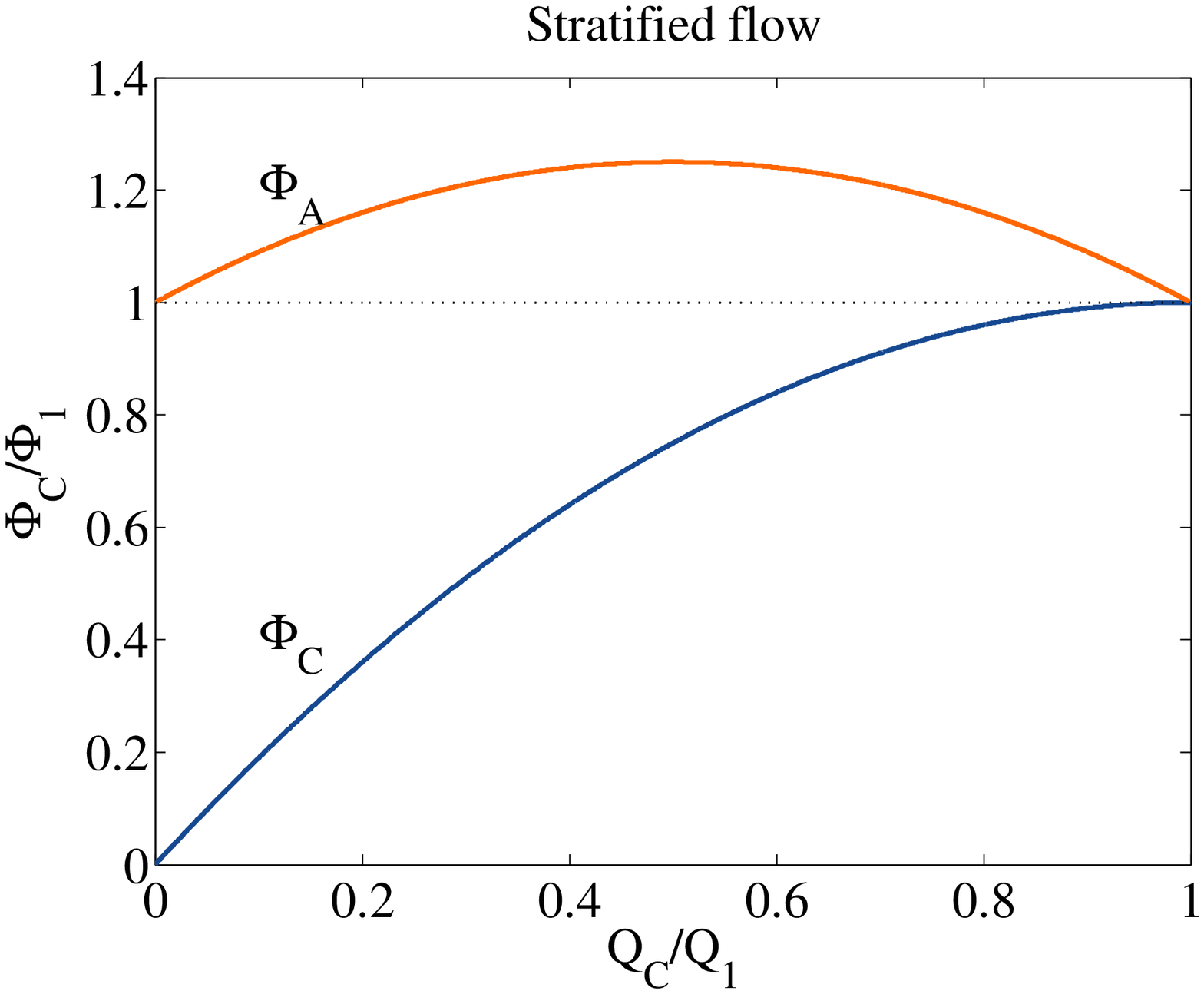}
        \end{subfigure}
\caption{Two examples of phase separation functions at a single node. Notation is for the upper node in Figure \ref{fig:schematic} when $Q_C > 0$. The normalized volume fraction in vessels $A$ and $C$  of the node is plotted as a function of the flow in vessel $C$ normalized by the inlet flow. The dotted line denotes the case with no phase separation. a) Empirical function for microvascular blood flow as defined by Equation \ref{eq:plasma_skim} and b) empirical function for stratified laminar flow as defined by Equation \ref{eq:baby_karst}.
}
\label{fig:visc_plasma}
\end{figure}

The phase separation rule for each node is a complex function which depends sensitively on the fluid system, the node geometry, and the inlet flow rate. The phase separation rule relates the downstream volume fractions in two daughter branches to the inlet flow state. For this work we explore the consequences of two different separation functions which are  valid for 1) microvascular blood flow and 2) stratified laminar flow of two Newtonian fluids. For microvascular blood flow, numerous authors have demonstrated this separation of red blood cells from plasma (\emph{i.e.}, plasma skimming) and developed simple empirical models to describe the effect \cite{Bugliarello:1964aa,Chien:1985aa,Dellimore:1983aa,Fenton:1985aa,Klitzman:1982aa,Pries:1989aa}. These empirical relations become part of the network model. In our previous work on networks with stratified laminar flow where the fluids remain as distinct phases, we measured the separation function for this system, demonstrated that we could compute the functions via 3D Navier-Stokes simulations, and developed an approximate one-parameter model for use in network modeling~\cite{karst2013}. In that work gravity was normal to the plane of the network flow. It has been shown that, unlike the effective viscosity model, the exact form of the separation function has a dramatic effect on the types of equilibrium and dynamic behavior that may be observed~\cite{Geddes:2007,karst2013}. The two sample empirical separation functions we use in this work are show in Figure \ref{fig:visc_plasma}.

It is important to realize that phase separation at a single network node is a common phenomena in many two fluid systems and that many other types of behaviors exist as noted in Section I. In network modeling, it is common to use  simple empirical functions with a single fit parameter which can be tuned to approximately model experimental data. It is recognized that such simple functions are limited in their accuracy, but they are useful in allowing for easy incorporation into analysis and providing some insight into expected experimental behavior. For example, a common fit function for microvascular flows is
\begin{equation}
\frac{\Phi_C }{\Phi_1 } =  \frac{ \left( \frac{Q_C}{Q_1}\right)^{p-1}}{ \left( \frac{Q_C}{Q_1}\right)^p + \left(1-\frac{Q_C}{Q_1}\right)^p} =
f \left( \frac{Q_C}{Q_1}\right),
\label{eq:plasma_skim}
\end{equation}
where $Q_C / Q_1$ is the normalized flow in branch $C$ shown in the schematic  of Figure \ref{fig:visc_plasma}.  We selected $p=2$ in Figure \ref{fig:visc_plasma}a; a typical value used in prior studies~\cite{Klitzman:1982aa}. For stratified flow a simple fit function which represents the basic form of the experimental data is,
\begin{equation}
\frac{\Phi_C}{\Phi_1} = 1- \gamma \left(1-\frac{Q_C}{Q_1}\right)^2 = f\left( \frac{Q_C}{Q_1}\right),
\label{eq:baby_karst}
\end{equation}
where in Figure \ref{fig:visc_plasma}b we selected $\gamma=1$,  which is observed in typical experimental data~\cite{karst2013}. In both cases the fit parameters $p$ and $\gamma$ depend on many of the other physical parameters in the system. We use the generic function $f$ to represent the phase separation constitutive law, whatever the physical system. For any function $f$, the volume fraction in vessel $A$ is connected to the function $f$ through conservation,  $Q_1 \Phi_1 = Q_A \Phi_A + Q_C \Phi_C $, which can be expressed as,
\begin{equation}
\frac{\Phi_A }{\Phi_1} = \frac{1  - f\left(\frac{Q_C}{Q_1}\right)\frac{Q_C}{Q_1}}{1-\frac{Q_C}{Q_1}}.
\end{equation}
A few remarks are worth making about the phase separation functions shown in Figure \ref{fig:visc_plasma}. Note that the phase separation function for microvascular blood flow is symmetric under the exchange $Q_1 \longleftrightarrow 1-Q_1$, \emph{i.e.}, it does not matter how we arrange the downstream vessels. The same is not true for the phase separation function of stratified flow; different arrangements of the downstream vessels results in different phase separation. In both cases we note that the volume fraction entering vessel $C$ is zero when $Q_C=0$, but this condition does not hold in any general sense. For all  phase separation functions $f(1)=1$ must hold.

\subsection{Governing equations}

We now develop a general network model based on conservation laws. We treat the viscosity function and phase separation function as constitutive laws which we must select in order to make concrete calculations of a real physical system. For all cases we use the viscosity law for mixed Newtonian fluids, Equation \ref{eq:viscosity}. In our network model we assume that the function $f$ is known by some means, either experiments or computational fluid dynamics. While we confine our results to Equations \ref{eq:plasma_skim} and \ref{eq:baby_karst}, the methods we develop  are general and can be applied to any $C^1$-smooth constitutive law for the physical system of interest.

We assume the volume fraction $\Phi_i(x_i,t)$ in vessel $i = A,B,C$ is governed by the first order wave equation
\begin{equation}
{\partial \Phi_i \over \partial t} + v_i {\partial \Phi_i \over \partial x_i} = 0, \quad 0 \leq x_i \leq \ell_i, \quad 0 \leq {t},
\label{eq:wave_dimensions}
\end{equation}
where $\ell_i$ is length of vessel $i$. The propagation velocity $v_i(t)$ in vessel $i$ is proportional to the volumetric flow rate $Q_i(t)$ in the vessel,
\begin{align}
v_i(t ) = {4 Q_i(t) \over \pi d_i^2},
\end{align}
where $d_i$ is the diameter of the vessel. At each node in the network the inlet flow rates equal the outlet flow rates, namely $Q_C = Q_1 - Q_A$ and $Q_B = Q_1 + Q_2 - Q_A$ where positive $Q_C$ is assumed to go from inlet $1$ to $2$. The flow rates may vary in time, but each is constant throughout the vessel. In this work we consider steady boundary conditions, thus $Q_1$, $\Phi_1$, $Q_2$, and $\Phi_2$ are constants.

To solve Equation \ref{eq:wave_dimensions} we need boundary conditions at the entrance of the three vessels. The boundary conditions are supplied by the conservation of each constituent at the node, namely,
\begin{align}
\Phi_A(0,t) &= \frac{\Phi_1 Q_1 - \Phi_C(0,t) Q_C(t) }{Q_1-Q_C(t)} \label{eq:BC1} \\ \nonumber \\
\Phi_B(0,t) & = \frac{\Phi_2 Q_2 + \Phi_C(\ell_c,t) Q_C(t) }{Q_2 + Q_C(t)} \label{eq:BC2}.
\end{align}
The third required boundary condition  depends upon the direction of $Q_C$. When the flow is such that $Q_C$ is positive, the boundary condition for vessel $C$  is given by $\Phi_C(0,t)=\Phi_1 f(Q_C/Q_1)$; see Figure \ref{fig:visc_plasma}b. When the flow is such that $Q_C$ is negative, the boundary condition in vessel $C$ is $\Phi_C(\ell_C,t)=\Phi_2 f(-Q_C/Q_2)$. Once the direction is established and the inlet volume fraction to vessel $C$ is determined by the phase separation function, Equations \ref{eq:BC1} and  \ref{eq:BC2} provide the inlet volume fractions to vessels $A$ and $B$.

The pressure drop across any vessel is given as  $\Delta P_i  = Q_i R_i$. In laminar flow, the hydraulic resistance of branch $i$, $R_i$, is a function of the spatially averaged  viscosity,
\begin{equation}
\bar{\mu}_i(t) = \frac{1}{\ell_i} \int_0^{\ell_i} \mu(\Phi_i(x,t))dx,
\end{equation}
through Poiseuille's law,
\begin{align}
R_i(t) = {128 \ell_i  \bar{\mu}_i(t) \over \pi d_i^4}.
\end{align}
Kirchoff's potential law applied around the network loop, \emph{i.e.}, $\Delta P_A = \Delta P_B + \Delta P_C$, provides an equation for the flow in $C$,
\begin{equation}
Q_C(t) = {Q_1R_A(t) - Q_2 R_B(t) \over R_A(t) + R_B(t) + R_C(t)}.\label{eqn:flow}
\end{equation}
The above formulation is a closed problem for the 1D wave propagation of volume fraction in the connected vessels of our network. It is worth noting that the model is symmetric under the exchange $Q_C \longleftrightarrow -Q_C$, $Q_1 \longleftrightarrow Q_2$, $\Phi_1 \longleftrightarrow \Phi_2$, and (vessel $A$) $\longleftrightarrow$ (vessel $B$).

\subsection{Dimensionless formulation}

A dimensionless version of the governing equations can be derived by scaling space and time according to
\begin{equation}
\hat{x}_i = {{x}_i \over \ell_i}, ~~~~~
\hat{t} = {Q_1 + Q_2 \over V_A + V_B + V_C} {t}, ~~~~~
\hat{Q}_i  = \frac{Q_i}{Q_1 + Q_2},
\end{equation}
so that each vessel's spatial dimension is normalized to its length,  time is scaled by the ratio of the total volumetric flow rate in the network to the total volume $V = V_A + V_B + V_C$ of the network, and flow rates are normalized to the total flow. The dimensionless governing equation for vessel $i$ is
\begin{align}
{\partial \Phi_i \over \partial \hat{t}} + \hat{Q}_i(t) {V \over V_i} {\partial \Phi_i \over \partial \hat{x}_i} = 0, \quad 0 \leq \hat{x}_i \leq 1, \quad 0 \leq \hat{t}. \label{eqn:pde}
\end{align}
The boundary conditions become,
\begin{align}
\Phi_A(0,t) &= \frac{\Phi_1 \hat{Q}_1     - \Phi_C(0,t) \hat{Q}_C(t) }{\hat{Q}_1 - \hat{Q}_C(t)}, \label{eq:BC}\\
\Phi_B(0,t) &= \frac{\Phi_2 (1-\hat{Q}_1) + \Phi_C(1,t) \hat{Q}_C(t) }{1 - \hat{Q}_1 + \hat{Q}_C(t)},
\end{align}
with the phase separation function at the appropriate node providing the final third boundary condition,
\begin{align}
\Phi_C(0,t)  &= \Phi_1 f\left(\frac{\hat{Q}_C}{\hat{Q}_1}\right)  ~~~~~~~~~ \mathrm{when} ~~~\hat{Q}_C>0, \\
\Phi_C(1,t)  &= \Phi_2 f\left(\frac{-\hat{Q}_C}{1-\hat{Q}_1}\right) ~~~~~ \mathrm{when} ~~~ \hat{Q}_C<0.
\label{eq:plasmaBC}
\end{align}
In dimensionless terms, the flow equation becomes
\begin{equation}
\hat{Q}_C = { \hat{Q}_1 r_A \bar{\mu}_A   - (1-\hat{Q}_1) r_B \bar{\mu}_B \over r_A \bar{\mu}_A +  r_B \bar{\mu}_B + r_C \bar{\mu}_C},
\label{eqn:flowEqn}
\end{equation}
where $r_i = 128 \mu_{\alpha} \ell_i / \pi d_i^4$ is the nominal resistance in vessel $i$, and $\bar{\mu}_i$ is the average relative viscosity in vessel $i$ as defined by Equation \ref{eq:viscosity}.

There are 8 dimensionless parameters that enter the problem. The network geometry introduces four parameters. Two of these are defined by the ratio of the nominal resistances, $r_A/r_C$ and $r_B/r_C$. The other two are defined by the ratio of the volume of the vessels, $V_A/V_C$ and $V_B/V_C$. In addition, there are three inlet parameters we are free to control, $\hat{Q}_1$, $\Phi_1$, and $\Phi_2$. The fluid system chosen determines the viscosity function, and the contrast between the two phases, $\mu_\beta/\mu_\alpha$, supplies another parameter. Finally, the phase separation function $f$ is critical to the behavior, though the function is set by the physical system and is not something we can easily control in a given physical experiment. The parameter space is quite large, thus direct numerical solution of the problem is not practical for spanning parameter space and motivates us to find a reliable method for tracking regions of stability and instability.

In what follows we use the dimensionless formulation, and for convenience we drop the ``hat'' notation.

\section{Equilibria}

At equilibrium, the volume fraction $\Phi_i(x_i,t)$ in branch $i$ is constant throughout the branch and equal to the entrance volume fraction, $\Phi_i(0,t)$. Equations \ref{eq:BC}--\ref{eqn:flowEqn} are sufficient to solve for the equilibrium flows and volume fractions. The viscosity functions and in turn hydraulic resistances can each be written as functions of the equilibrium flow rate $Q_C$. Equation \ref{eqn:flowEqn} therefore defines a nonlinear equation in $Q_C$, and multiple solutions are possible. The parameter space is still large and consists of $r_A/r_C$, $r_A/r_B$, $Q_1$, $\Phi_1$, $\Phi_2$, and $\mu_{\beta}/\mu_{\alpha}$.

We have explored the equilibrium problem for the 3-node network in prior publications. Gardner \emph{et al.}~\cite{Gardner:2010} theoretically studied the problem in the context of microvascular blood flow and demonstrated that multiple equilibrium states were indeed possible. The observation that multiple equilibria were possible in regimes where no phase separation takes place motivated us to design a table-top experiment using water and sucrose solution \cite{geddes2010a}. In this work we derived a simple condition for the onset of multiple equilibria, and confirmed the predictions in the laboratory \cite{geddes2010a}. More recently, Karst \emph{et al.}~\cite{karst2013} designed an experiment using fluids undergoing laminar stratified flow to attempt to mimic the phase separation effect in microvascular blood flow. In that work we predicted and  observed multiple equilibria and derived a simple criteria for its onset. In this current paper, we focus on the problem of stability and dynamics. However, for completeness we synthesize prior results on equilibria here using our current models and terminology. We refer the interested reader to the above publications for more details.

\begin{figure}[]
    \begin{subfigure}[]
                \centering
                \includegraphics[width=66.5mm]{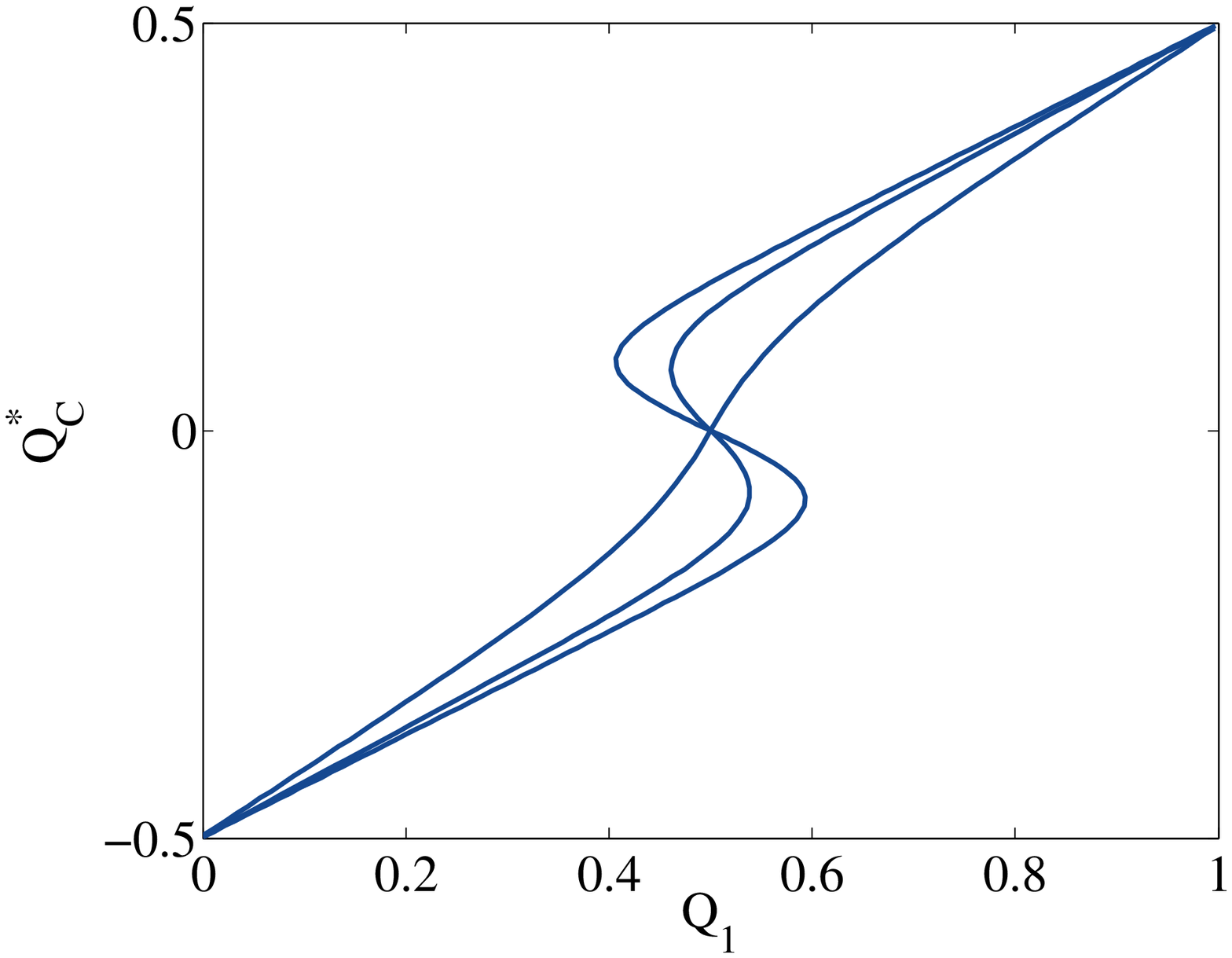}
        \end{subfigure}
        \begin{subfigure}[]
                \centering
                \includegraphics[width=65mm]{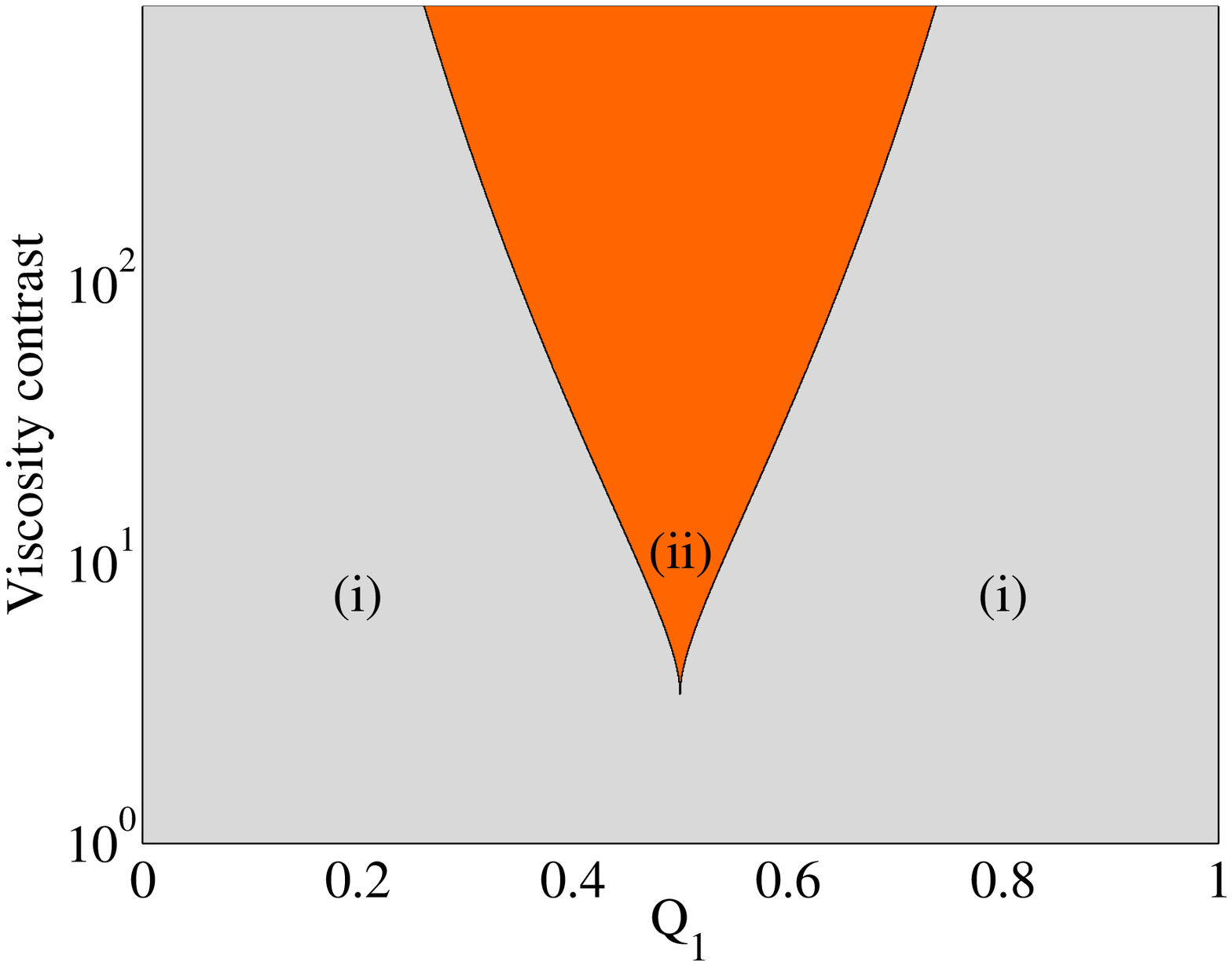}
        \end{subfigure}
\caption{(a) Equilibrium curves for viscosity contrast of 2, 10, and 30 in the microvascular blood flow model. The equilibrium curve is single-valued when the viscosity contrast is 2. Multiple equilibria are created via a saddle-node bifurcation at $(Q_1,Q_C^*) = (0.5,0)$. When the viscosity contrast is set to 10, we observe a small window of multiple equilibria about $Q_1 = 0.5$, and as the viscosity contrast is increased to 30, this window widens. (b) Phase diagram in the $Q_1 \times (\mu_\beta / \mu_\alpha)$ parameter space. In the gray region (i), there exists a single equilibrium. In the orange region (ii), there exist multiple equilibria. Regions (i) and (ii) are delineated by the saddle-node bifurcation curve (black) which emerges from (0.5,3.4). }
\label{fig:mvEq}
\end{figure}

In Figure \ref{fig:mvEq}a we show three sample equilibrium curves in the $(Q_1,Q_C^*)$ plane. We have chosen $r_A/r_C = 4(2.5)^4/3 \approx 52.1$, $r_A/r_B = 1$, $\Phi_1 = \Phi_2 = 0.82$, and the three curves correspond to viscosity contrast $\mu_{\beta}/\mu_{\alpha} = 2, 10,$ and $30$. Here we are using Equation \ref{eq:plasma_skim} for the phase separation function $f$ from microvascular blood flow. The parameters selected here are relevant later in our analysis of the dynamics. There exist multiple equilibria if for any given value of $Q_1$ there exist multiple values of $Q_C^*$. For the chosen parameter values, there is a single equilibrium for a viscosity contrast of $2$, but multiple equilibria for contrasts of $10$ and $30$. The window of multiple equilibria grows with increasing viscosity contrast.

While the width of the multiple equilibria window involves an in-depth calculation, the onset point is relatively straight-forward to calculate. Notice from Figure \ref{fig:mvEq}a that multiple equilibria are born in a saddle-node bifurcation when the equilibrium curve folds over at $Q_C^*=0$. A condition for onset can therefore be obtained by setting $dQ_1/dQ_C^* = 0$, which yields,
\begin{align}
 R_A + R_B + R_C = \ln\left({\mu_{\beta} \over \mu_{\alpha}}\right) \left(R_A (\Phi_1 - \Phi_C) + R_B (\Phi_2 - \Phi_C)\right).
\end{align}
Recall that $R_A= r_A \bar{\mu}_A$, where $r_A$ depends only upon the geometry ($d_A$ and $\ell_A$) of  vessel A while $R_A$ depends upon the phase distribution within the network. Thus to evaluate the hydraulic resistances $R_i$ we must know the network geometry and the phase distribution inside the network when $Q_C^*=0$. If the inlets are not equal fluids we must be careful to consider the above criteria as  $Q_C^* \rightarrow 0^{+}$ and $Q_C^* \rightarrow 0^{-}$.

We can simplify the multiple equilibrium criteria for the two cases considered in this paper. First, we limit our study to cases where we drive the network with identical inlet fluids, $\Phi_1 = \Phi_2$, and we do not need to consider the direction with which we approach $Q_C^* \rightarrow 0$. Second, for the two-phase separation functions  examined in this paper, the volume fraction in vessel $C$ is zero when $Q_C=0$; $f(0) = 0$ in both empirical phase separation functions given by Equations \ref{eq:plasma_skim}  and \ref{eq:baby_karst}. In this particular case, the condition for multiple equilibrium becomes
\begin{align}
 R_A + R_B+R_C = \ln\left({\mu_{\beta} \over \mu_{\alpha}}\right) \left( R_A  + R_B  \right) \Phi_1.
\end{align}
Since at $Q_C^*=0$, $\Phi_A=\Phi_B=\Phi_1$ and $\Phi_C=0$, the criteria can be further reduced to
\begin{align}
\left( 1 + \frac{1}{\mu_1} \frac{r_C}{ r_A + r_B}  \right) =  \ln\left( \mu_{\beta} \over \mu_{\alpha}\right) \Phi_1 = \ln\left(\mu_1 \right),
\end{align}
where $\mu_1$ is the relative viscosity of the inlet fluid. For the network geometry used in Figure \ref{fig:mvEq}a, $r_C \ll r_A + r_B$, thus the criteria for multiple equilibria approximately reduces to $\ln{\mu_1} = 1$, or  $\mu_{\beta}/ \mu_{\alpha} =\mathrm{e}^{1/\Phi_1}\approx 3.4$.

In Figure \ref{fig:mvEq}b we show the region of multiple equilibria in the $Q_1 \times (\mu_{\beta}/\mu_{\alpha})$ plane for the parameters previously discussed. Notice that the onset point agrees with the above calculation and occurs at  $(0.5,3.4)$. As the viscosity contrast is increased the width of the window increases. Changing the network geometry and inlet fluids changes the details of the multiple equilibria window but not its existence. In the rest of this paper, we consider the stability of the equilibrium solutions and the resulting nonlinear dynamics.

\section{Linearization and the characteristic equation}
We assume that the network is initially in equilibrium, \emph{i.e.}, $Q_i(t) = Q_i^*$, $\Phi_i(x_i,t) = \Phi_i^*$ for all $t < 0$ and $i = A,B,C$. We introduce perturbations  beginning at time $t = 0$ on the flow rates $Q_i(t)$ and volume fraction profiles $\Phi_i(x_i,t)$ so that
\begin{align}
\Phi_i(x_i,t) &= \Phi_i^*(1 + \Delta \Phi_i(x_i,t)) \label{eqn:perturbPhi} \\
Q_i(t) &= Q_i^*(1 + \Delta Q_i(t)). \label{eqn:perturbQ}
\end{align}
Substituting Equations \ref{eqn:perturbPhi} and \ref{eqn:perturbQ} into the appropriate governing equations and retaining only the linear terms results in a first order wave equation describing the propagation of the volume fraction perturbation in each branch,
\begin{align}
{\partial \over \partial t}  \Delta \Phi_i + {1 \over \tau_i} {\partial \over \partial x_i}  \Delta \Phi_i = 0,
\label{eq:wave}
\end{align}
where $\tau_i = Q_i^* V / V_i$ is the dimensionless steady state transit time in branch $i$. An expression for the flow perturbation can be computed by expanding the flow equation about the equilibrium,
\begin{align}
\Delta Q_C(t) &= {Q_A^* R_A^* \Delta R_A(t) - Q_B^* R_B^* \Delta R_B(t) - Q_C^* R_C^* \Delta R_C(t) \over Q_C^* \sum_i R_i^*}.
\label{eq:flow}
\end{align}
Relative perturbations to the resistance in each branch is determined by
\begin{align}
R_i(t) &= {128 \ell_i \over \pi d_i^4} \int_0^1 \mu_i\left(\Phi_i^*(1 + \Delta \Phi_i(x_i,t))\right) dx_i \nonumber \\
&= {128 \ell_i \over \pi d_i^4} \int_0^1 \mu_i(\Phi_i^*) + \left.{d \mu_i \over d \Phi_i} \right|_* \Phi_i^* \Delta \Phi_i(x_i,t) dx_i \nonumber \\
&= R_i^* + R_i^*  \Phi_i^* \left.{d \ln(\mu_i) \over d \Phi_i} \right|_* \int_0^1 \Delta \Phi_i(x_i,t) dx_i \nonumber \\
\Rightarrow \Delta R_i (t) &= \Phi_i^* \left.{d \ln(\mu_i) \over d \Phi_i} \right|_* \int_0^1 \Delta \Phi_i(x_i,t) dx_i.
\label{eq:resist}
\end{align}
Finally, perturbations to the boundary conditions are required. Without loss of generality we assume that the flow in $C$ is from inlet 1 to inlet 2 and the perturbation to the volume fraction entering $C$ is then
\begin{align}
\Phi_C(0,t) &= \Phi_1 f\left({Q_C^*(1 + \Delta Q_C(t)) \over Q_1}\right) \nonumber \\
&= \Phi_C^* + \left.\left({\Phi_1 Q_C \over Q_1} f'\right)\right|_* \Delta Q_C(t) \nonumber \\
\Rightarrow \Delta \Phi_C(0,t) &=  \left.\left({Q_C \over Q_1} {f' \over f}\right)\right|_* \Delta Q_C(t),
\label{eq:vesselC}
\end{align}
where $f'$ is the derivative of the plasma skimming function $f$. The perturbations to the volume fraction entering $A$ and $B$ are given by mass fraction. For vessel $A$ we have
\begin{align}
Q_A^* \Phi_A^* \Delta \Phi_A(0,t) &= Q_C^* (\Phi_A^* - \Phi_C^*) \Delta Q_C(t)  - Q_C^* \Phi_C^* \Delta \Phi_C (0,t),
\label{eq:vesselA}
\end{align}
and for vessel $B$ we have
\begin{align}
Q_B^* \Phi_B^* \Delta \Phi_B(0,t) &= Q_C^* (\Phi_C^* - \Phi_B^*) \Delta Q_C(t)  + Q_C^* \Phi_C^* \Delta \Phi_C (1,t).
\label{eq:vesselB}
\end{align}
Equations \ref{eq:wave} -- \ref{eq:vesselB} constitute the linearized equations. We assume traveling wave solutions of the form
\begin{align}
\Delta Q_C(t) &= \Delta q_C e^{\lambda t} \\
\Delta \Phi_i(x_i,t) &= \Delta \phi_i e^{\lambda(t - \tau_i x_i)},
\end{align}
which automatically satisfy Equation \ref{eq:wave}. Substituting into Equation \ref{eq:resist} and integrating gives
\begin{align}
\Delta R_i(t) = \Delta \phi_i  u_i e^{\lambda t}  .
\end{align}
where
\begin{align}
u_i = \Phi_i^* \left.{d \ln(\mu_i) \over d \Phi_i} \right|_* {(1- e^{-\lambda \tau_i}) \over \lambda \tau_i}.
\end{align}
Further substitution into Equation \ref{eq:flow} results in
\begin{align}
\Delta q_C =  {R_A^* Q_A^* u_A \Delta \phi_A - R_B^* Q_B^* u_B \Delta \phi_B  -  R_C^* Q_C^* u_C \Delta \phi_C \over Q_C^* \sum_i R_i^*}.
\label{eq:linear1}
\end{align}
Substituting into Equation \ref{eq:vesselC} results in
\begin{align}
\Delta \phi_C = \left.\left({Q_C \over Q_1} {f' \over f}\right)\right|_* \Delta q_C.
\end{align}
Finally, substitution into Equations \ref{eq:vesselA} and \ref{eq:vesselB} gives
\begin{align}
Q_A^* \Phi_A^* \Delta \phi_A =  Q_C^* (\Phi_A^* - \Phi_C^*) \Delta q_C - Q_C^* \Phi_C^* \Delta \phi_C
\end{align}
and
\begin{align}
Q_B^* \Phi_B^* \Delta \phi_B =  Q_C^* (\Phi_C^* - \Phi_B^*) \Delta q_C +Q_C^* \Phi_C^* \Delta \phi_C e^{-\lambda \tau_C}.
\label{eq:linear4}
\end{align}
Equations (\ref{eq:linear1})-(\ref{eq:linear4}) constitute 4 linear equations in the 4 unknowns $\Delta q_A, \Delta \phi_A, \Delta \phi_B,$ and $\Delta \phi_C$. Non-trivial solutions exist if and only if the following characteristic equation has roots,
\begin{align}
\chi(\lambda) = a \left({1 - e^{-\lambda \tau_A} \over \lambda \tau_A}\right) + (b + de^{-\lambda \tau_C})\left({1 - e^{-\lambda \tau_B} \over \lambda \tau_B}\right) + c\left({1 - e^{-\lambda \tau_C} \over \lambda \tau_C}\right) - 1,
\label{eqn:char}
\end{align}
where the coefficients are given by
\begin{align}
a &= -\left((\Phi_C-\Phi_A) + \Phi_C {Q_C \over Q_1} {f' \over f} \right) {R_A \over \sum_i R_i} {d \ln(\mu_A) \over d \Phi_A}, \\
c &= -\Phi_C {Q_C \over Q_1}  {f' \over f}  {R_C \over \sum_i R_i} {d \ln(\mu_C) \over d \Phi_C}, \\
b &= -(\Phi_C - \Phi_B)   {R_B \over \sum_i R_i} {d \ln(\mu_B) \over d \Phi_B}, \\
d &= -\Phi_C {Q_C \over Q_1} {f' \over f} {R_B \over \sum_i R_i} {d \ln(\mu_B) \over d \Phi_B}.
\end{align}
and we have dropped the * for convenience.

The characteristic equation has three delay times, but is composed of linear combinations of four transcendental functions. Two of these arise from the propagation delay in vessel $A$ and vessel $C$ (coefficients ``a" and ``c"). A third term arises due to perturbations in the flow entering vessel $B$ (coefficient ``b"). The last contribution arises due to perturbations to the volume fraction in vessel $C$ which propagate into and through vessel $B$ (coefficient ``d"). It is worth noting that in the absence of nonlinear viscosity, all of the coefficients are zero and the equilibrium is stable. Furthermore, no plasma skimming would imply $f'=0$ and $\Phi_C = \Phi_A$ so that only the ``b" coefficient would remain. It is straightforward to show that only real roots exist and thus oscillatory dynamics are ruled out. Nonlinear viscosity and plasma skimming are therefore necessary for the emergence of oscillatory behavior.

A root $\lambda = \sigma + i \omega$ of the characteristic equation satisfies the relations
\begin{align}
R(\sigma,\omega) = \Re(\chi(\sigma + i\omega)) &= 0 \label{eqn:eigRe}  \\
I(\sigma,\omega) = \Im(\chi(\sigma + i\omega)) &= 0.\label{eqn:eigIm}
\end{align}
We will see that these relations are useful for identifying Hopf bifurcations that can be used as starting points for numerical continuation through the large parameter space of the system.

\section{Results}

The network model includes 8 dimensionless parameters as well as the constitutive laws for viscosity and phase separation. It is difficult to make general predictions without selecting a set of constitutive laws since these relations critically determine the system behavior. Since general statements about any arbitrary system are difficult to make, we present two physically realistic systems to demonstrate the methodology for analyzing stability. In Example 1 we take the well-studied problem of microvascular blood flow, and in Example 2 we take stratified laminar flow of two Newtonian fluids, a system for which we have conducted prior equilibrium experiments  \cite{karst2013}.

\subsection{Example 1: Microvascular blood flow}

For our first example,  we use the phase separation model for microvascular flow, Equation \ref{eq:plasma_skim} with $p = 2$. We use the simple Arrhenius law for viscosity in the vessels after the initial splitting at the inlets, Equation \ref{eq:viscosity}. The Arrhenius law has the basic functional form as the empirical laws for blood viscosity \cite{Geddes:2007}. For the network we use parameters $\Phi_1 = \Phi_2 = 0.82$;  $d_A = d_B = 1, d_C = 2.5$; $\ell_A = \ell_B = 1, \ell_C = 0.75$ unless otherwise noted. In dimensionless terms $r_A/r_C = 52.1$, $r_A/r_B=1$, $V_A/V_C=0.213$, and $V_A/V_B=1$.

\begin{figure}[!t]
\centering
\includegraphics[width=100mm]{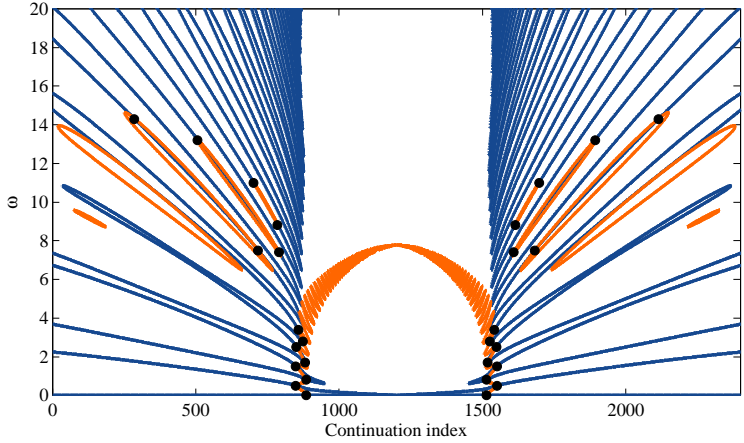}
\caption{Zero contours of Equation \ref{eqn:eigRe} (blue) and Equation \ref{eqn:eigIm} (orange) with $\sigma = 0$ and  a viscosity contrast of $50$ in the microvascular blood flow model. Each intersection (black dot) indicates a Hopf bifurcation of frequency $\omega$ occurs at the $(Q_C^*,Q_1)$ pair associated with the continuation index.
}\label{fig:hopfContours}
\end{figure}

The traditional approach to detect Hopf bifurcations is to monitor the test function defined by the product of the imaginary components of the eigenvalues along the continuation of an equilibrium. In the systems with transcendental characteristic equations in which the eigenvalues can not be directly computed, a more powerful tool can be applied by monitoring the test function $\det(2 J \odot I_n)$, where $J$ is the Jacobian of the equilibrium relation and $\odot$ denotes the bialternate product. Here, we employ a more specialized approach in order to simultaneously detect a Hopf bifurcation and determine its frequency.

At equilibrium, the hydraulic resistances in each branch are functions of the equilibrium flow rate $Q_C^*$. We can therefore rewrite Equation \ref{eqn:flow} as $Q_C^* = \psi(Q_C^*)$. To track an equilibrium through parameter space, we parameterize $\psi$ by $Q_1$ and perform numerical continuation on the equilibrium relation
\begin{align}
F_E(Q_C^*,Q_1) = \psi(Q_C^*,Q_1) - Q_C^*,
\end{align}
forming a parametric equilibrium curve $\beta(s)$ in the $Q_1 \times Q_C^*$ plane. Hopf bifurcations can be identified along the equilibrium curve by monitoring the relation defined by substituting $\sigma = 0$ in Equations \ref{eqn:eigRe} and \ref{eqn:eigIm}. Since the values of  $a,b,c,d,$ and the steady state transit times $\tau_i$ are fixed at each $(Q_1,Q_C^*)$ pair along the continuation $\beta(s)$, we can define $R(s,\omega)$ and $I(s,\omega)$ to be the left sides of Equations \ref{eqn:eigRe} and \ref{eqn:eigIm} with $\sigma = 0$, respectively, without loss of generality. Then any intersection of the zero contours of $R(s,\omega)$ and $I(s,\omega)$ indicates a Hopf bifurcation of frequency $\omega$ occurs at the $(Q_1,Q_C^*)$ pair associated with index $s$. We  see an implementation of this methodology with a viscosity contrast of 50 in Figure \ref{fig:hopfContours}. This figure is horizontally symmetric because the underlying network is geometrically symmetric.  We observe a collection of low frequency Hopf bifurcations that occur near the saddle-node bifurcations located at indices 886 and 1515 in Figure \ref{fig:hopfContours}. We also observe pairs of higher frequency Hopf bifurcations that occur away from the saddle-node bifurcations. As the viscosity contrast is increased, additional bands of instability appear, and these bands grow to encompass the entirety of the upper and lower branches of the equilibrium curve.

\begin{figure}
    \begin{subfigure}[]
                \centering
                \includegraphics[width=65mm]{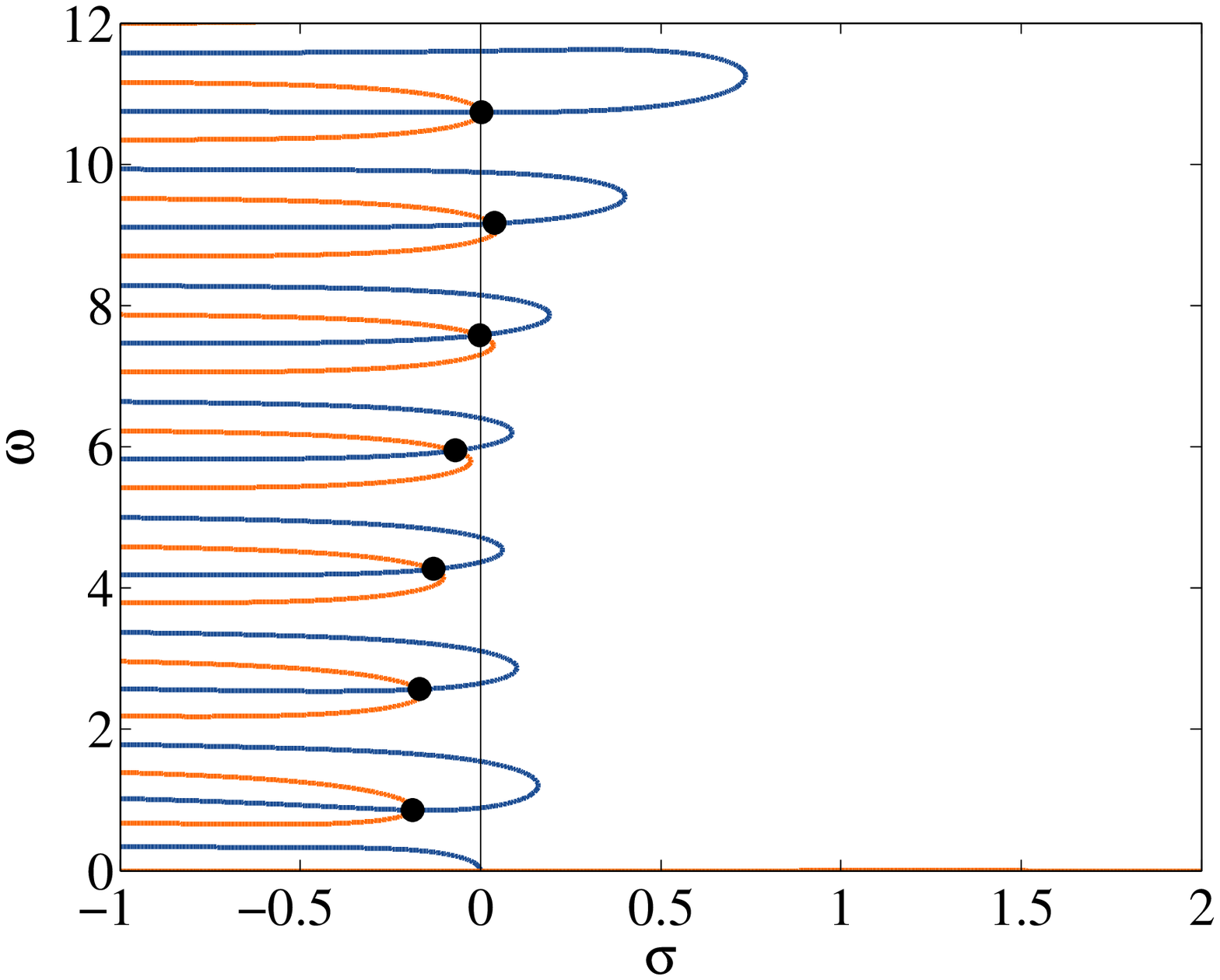}
        \end{subfigure}
        \begin{subfigure}[]
                \centering
                \includegraphics[width=68.5mm]{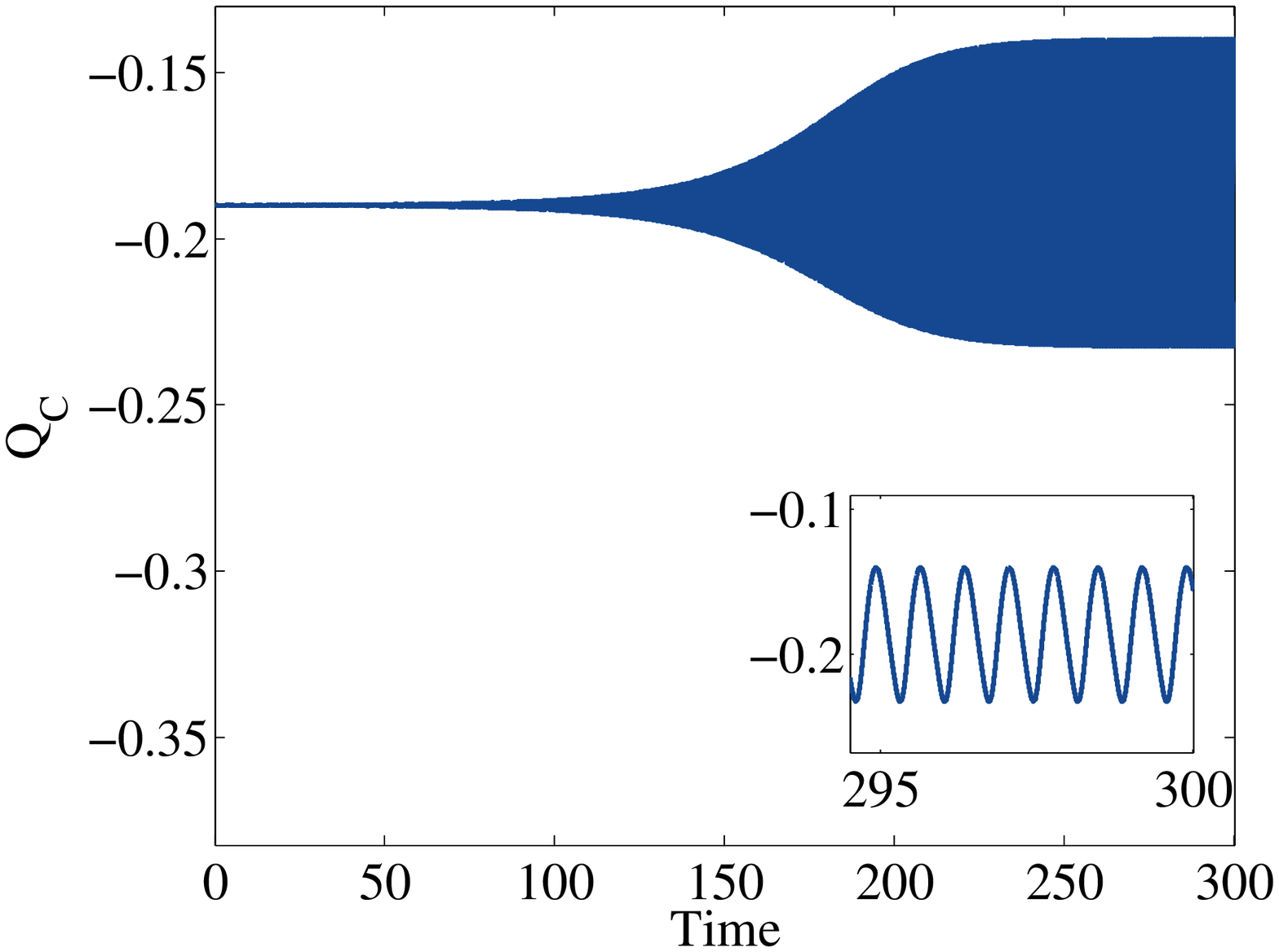}
        \end{subfigure}
        \caption{
        (a) Associated zero contours of Equation \ref{eqn:eigRe} (blue) and Equation \ref{eqn:eigIm} (orange) in the microvascular blood flow model. Each intersection (black dot) indicates a solution $\lambda = \sigma + i\omega$ to the characteristic equation \ref{eqn:char}. Note the sole eigenvalue with positive real part is $\lambda \approx 0.04 +  9.16i$. (b) Limit cycle about equilibrium $(Q_1,Q_C^*) = (0.5,-0.19)$ computed from direct simulation. The period of the oscillation agrees with the analysis.}\label{fig:eig}
 \end{figure}

We can confirm that Figure \ref{fig:hopfContours} accurately predicts the presence of sustained oscillations through direct numerical simulation. As an example, we choose the equilibrium pair $(Q_1,Q_C^*) = (0.5,-0.19)$ which is located in the left-most band of instability in Figure \ref{fig:hopfContours}. The eigenvalue-based prediction is shown in  Figure \ref{fig:eig}a. Here we plot the zero contours of Equations \ref{eqn:eigRe} and \ref{eqn:eigIm} in the $\sigma \times \omega$ plane so that an intersection of the contours at some $(\sigma,\omega)$ pair indicates that $\lambda = \sigma + i \omega$ is a solution to the Equation \ref{eqn:char}. Note the dominant eigenvalue has positive real part and imaginery part $\omega \approx 9.16$. We then perform a direct numerical simulation of Equation \ref{eqn:pde} (with appropriate boundary conditions) at the same parameters. We initialize the simulation to the equilibrium state and provide a small numerical perturbation. A  limit cycle grows from the unstable equilibrium solution as seen in Figure \ref{fig:eig}b.  When the system reaches a periodic steady state the limit cycle  has period ${T} \approx 0.717$, which corresponds to a dimensionless angular frequency $\omega = 2\pi / {T}\approx 8.84$, in good agreement with our linear prediction.
If we check the frequency in the simulation earlier when the amplitude is infinitesimal the frequency matches the linear analysis exactly.
We also confirm that our predicted growth rate matches the simulation.

\begin{figure}[!h]
\centering
\includegraphics[width=90mm]{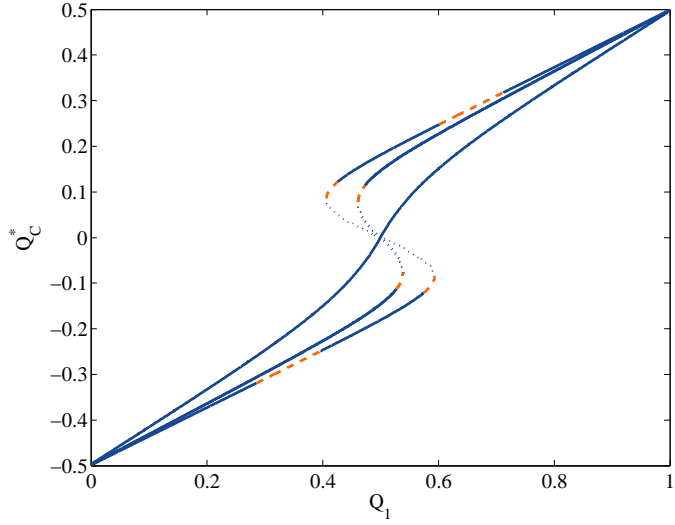}
\caption{Equilibrium curves of the microvascular blood flow model for viscosity contrasts of 2, 10, and 30. Solid regions of the equilibrium curve represent stable equilibria, while dotted region represent unstable equilibria. Dashed regions indicate the existence of a limit cycle.}\label{fig:hopfHighlight}
\end{figure}

We can begin to form intuition about the presence and location of Hopf bifurcations by varying the viscosity contrast for a fixed geometry and tracking the associated bands of instability. An example  is shown  in Figure \ref{fig:hopfHighlight}. Here we plot the equilibrium curve in the $Q_1 \times Q_C^*$ plane at three values of the viscosity contrast. This is the same figure and parameters as Figure \ref{fig:mvEq}a with the stability information superimposed. These curves are experimentally relevant as one can build a fixed network and then adjust the relative flow of the two inlets to move left and right along the $x$-axis \cite{karst2013}. Experimentally we can adjust the inlet fluids to adjust to viscosity, here the three curves represent the equilibrium solution for viscosity contrasts of 2, 10, and 30. When the viscosity contrast is 2, the equilibrium curves are single-valued and there are no Hopf bifurcations. At a viscosity contrast of 10, the equilibrium curve becomes multi-valued over a small range around $Q_1=0.5$. For  this range of $Q_1$  there are two possible states, one with positive and negative $Q_C^*$. We also see a region of instability emerges right at the location where the curves fold over. This Hopf bifurcation is at  low frequency and in numerical simulations we find that there is no stable limit cycle. The amplitude of oscillation grows until the system flips to the other  stable state on the equilibrium curve.

As we increase the viscosity contrast to 30 the region of multiple equilibrium grows and a new region of instability emerges along the equilibrium curve. This region is a  high frequency oscillation which results in a stable limit cycle as seen in Figure \ref{fig:eig}b. For the viscosity contrast of 30, the picture is that as we experimentally move continuously from $Q_1=0$ to $Q_1=1$ we would start by observing a single, stable, equilibrium flow state with negative $Q_C^*$. As we increase $Q_1$ we would see a limit cycle oscillation emerge around $Q_1\approx0.286$ which would persist until $Q_1\approx0.397$. Since this limit cycle exists outside the region of bistability, there is no other state for the system to move toward. After $Q_1$ is increased beyond 0.37 the limit cycle disappears and the system returns to a single stable equilibrium state with negative $Q_C^*$. At $Q_1\approx 0.575$ the large amplitude oscillation emerges and kicks the system to the other stable equilibrium state with positive $Q_C^*$.  As we continue to increase $Q_1$ the oscillations would emerge again at $Q_1\approx0.603$, this time with positive $Q_C^*$. Finally at $Q_1\approx0.714$ we would return to a single, stable equilibrium with positive $Q_C^*$. In this example the curves are symmetric about $Q_1=0.5$ because the network geometry is symmetric.

\begin{figure}[!h]
\centering
\includegraphics[width=120mm]{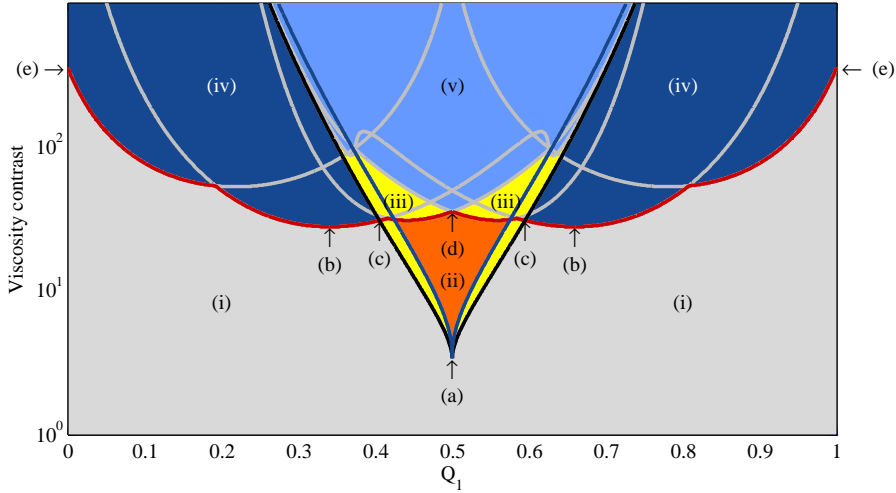}
\caption{Phase diagram in $Q_1 \times (\mu_\beta / \mu_\alpha)$
 parameter space for the microvascular blood flow model. In the gray region (i), the system exhibits a unique equilibrium state. In the orange region (ii), two stable equilibrium states exist.
The yellow region (iii) represents parameters which support one unstable oscillation and  one stable state. The dark blue region (iv) represents configurations in which have a single  oscillatory state. Networks in the light blue region (v) support two oscillatory states. The regions are separated by curves marking saddle-node bifurcations (black curves), the lowest frequency Hopf bifurcation (blue curve), and higher frequency Hopf bifurcations (red/gray curves). At high viscosity contrasts the instability is comprised of  multiple frequencies.}
  \label{fig:bifurc}
\end{figure}

The region of instability  changes as we increase the viscosity contrast. Generally, the window with multiple equilibrium states and the regions of instability increase with viscosity contrast. This behavior is demonstrated in the phase diagram of Figure \ref{fig:bifurc} which is an expansion of the phase diagram shown previously in Figure  \ref{fig:mvEq}b. Here we have identified both saddle-node and Hopf bifurcations in Figure \ref{fig:hopfContours} and used the relations
\begin{align*}
F_S(Q_1,Q_C^*,\mu,\omega) = \begin{bmatrix} \psi - Q_C^* \\ {d \psi \over dQ_C^*} - 1 \end{bmatrix}, \quad F_H(Q_1,Q_C^*,\mu,\omega) = \begin{bmatrix} \psi- Q_C^* \\ R(0,\omega) \\ I(0,\omega) \end{bmatrix},
\end{align*}
to track the saddle-node and Hopf bifurcations, respectively, through parameter space. Note that this phase diagram is symmetric about $Q_1=0.5$ due to the symmetry of the network geometry.

The Hopf and saddle-node continuation curves delineate several regions of behavior. If we start with a low  viscosity contrast, \emph{i.e.}, less than $3.4$,  we have single valued equilibrium curve for any $Q_1$. As we increase the viscosity contrast, multiple equilibrium behavior emerges from $Q_1=0.5$ (point a). As soon as multiple equilibrium exists, a small window of instability emerges right at the point that the equilibrium curves fold over. This is a narrow region of a low frequency, large amplitude oscillation which will generally kick the system to the stable part of the multiple equilibrium curve. Recall the behavior for a viscosity contrast of 10 from Figure \ref{fig:hopfHighlight}. As we increase the viscosity contrast to $27.8$, a small region of high frequency instability emerges at $Q_1\approx0.33$ and $Q_1\approx 1-0.33 = 0.67$ (point b). This first region of high frequency instability  represent the emergence of a single frequency limit cycle. The emergence of this instability occurs  outside the multiple equilibrium region, thus the system must oscillate around the equilibrium point. This  behavior was seen at a viscosity contrast of 30 in Figure \ref {fig:hopfHighlight}.

As we increase the viscosity contrast to $30.7$, the instability curve crosses into the region of multiple equilibria (point c). In this region, we find that the system may tend to the oscillatory solution with positive (negative) $Q_C^*$ or the stable solution with negative (positive) $Q_C^*$,  depending on the initial condition. As the viscosity contrast is increased to $35.5$, the Hopf bifurcation curves associated with the positive and negative $Q_C^*$ cross at $Q_1=0.5$. At this point we have the  co-existing limit cycles at $Q_1=0.5$ (point d); the system has two possible limit cycles one with positive and another with negative $Q_C^*$. We also see that at this viscosity value we have multiple frequency components meaning more complex dynamics are expected. Finally, as we increase the viscosity contrast to $352$, the region of instability reaches the boundary where $Q_1=0$ and $Q_1=1$ (point e); thus instability encompasses the entire range of inlet flow rates. All values of $Q_1$ are expected to be unstable and if we are inside the multiple equilibrium region we expect to always find co-existing limit cycles.

\begin{figure}[!h]
\includegraphics[width=125mm]{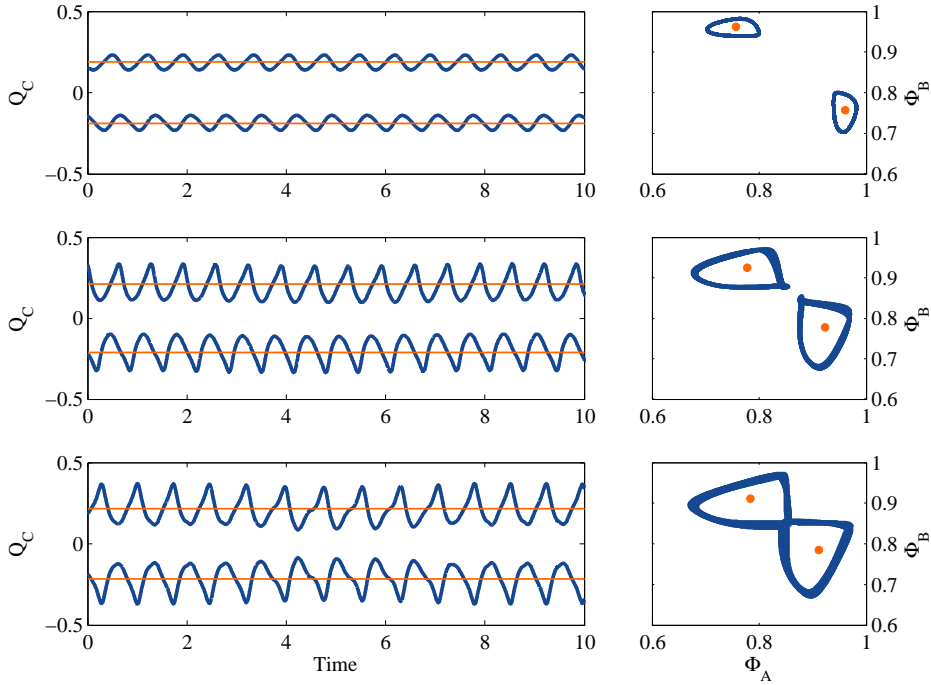}
\caption{Time series and associated phase plot from the direct numerical simulation of the microvascular blood flow model for different values of viscosity contrast, 50, 500, and 1500 from top to bottom. In the upper two plots we see coexisting limit cycles with increasing amplitude.  In each phase plot the equilibrium solution is shown as the dot. }\label{fig:increaseDelta}
\end{figure}

The presence of Hopf bifurcations is strongly dependent on the viscosity contrast $\mu_\beta/ \mu_\alpha$, and it not surprising that tuning this parameter also affects the amplitude and frequency of the associated oscillations. We saw in Figure \ref{fig:bifurc} that at high viscosity contrast we could have coexisting limit cycles and oscillations with multiple frequency components. In these cases our linear analysis can not tell us the complete dynamics so we use direct numerical simulation of Equation \ref{eqn:pde} to explore the final dynamics. In Figure \ref{fig:increaseDelta} we plot three time series of the flow in the middle branch $Q_C(t)$. As viscosity contrast  is increased from 50 to 500 to 1500, the amplitude of the limit cycle grows considerably.


\subsection{Example 2: Stratified laminar flow}

For example 2,  we use the phase separation model for stratified laminar flow, Equation \ref{eq:baby_karst} with $\gamma = 1$. We again use the simple Arrhenius law for viscosity in the vessels after the initial splitting at the inlet. This viscosity law could be realized in experiments if mixing was induced after the initial inlet split or if the tubes A, B, and C were long enough to allow molecular diffusion to mix the two phases. For the network we use similar parameters as the previous example $\Phi_1 = \Phi_2 = 0.8$;  $d_A = 1, d_B = 0.5, d_C = 2.5$; $\ell_A = \ell_B = 1, \ell_C = 0.75$ unless otherwise noted. In dimensionless terms, $r_A/r_C=52.1$, $r_A/r_B=\frac{1}{16}$, $V_A/V_C=0.213$ and $V_A/V_B=4$. Note that in comparison to the microvascular example we have broken the symmetry of the diameters in vessels $A$ and $B$. While we have no definite proof, we have been unable to detect any Hopf bifurcations in a symmetric network subject to stratified laminar flow.

\begin{figure}[!t]
\begin{center}
\includegraphics[width=100mm]{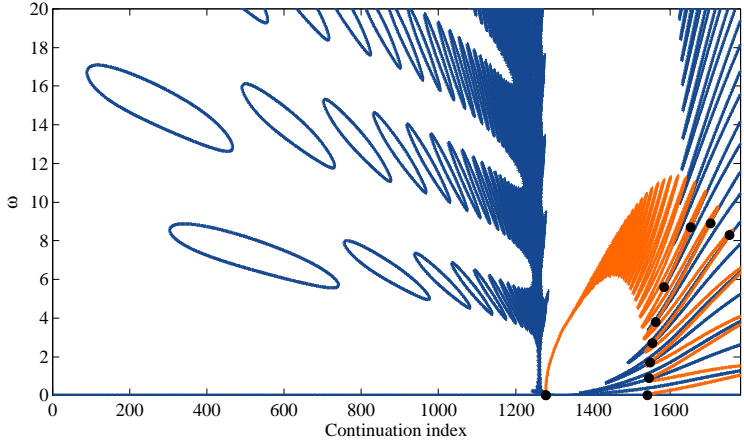}
\end{center}
\caption{Zero contours of Equation \ref{eqn:eigRe} (blue) and Equation \ref{eqn:eigIm} (orange) with $\sigma = 0$ and a viscosity contrast of 50 in the stratified laminar flow model. Each intersection (black dot) indicates a Hopf bifurcation of frequency $\omega$ occurs at the $(Q_1,Q_C^*)$ pair associated with the continuation index.}\label{fig:hopfContours_strat}
\end{figure}

We apply the same technique discussed in Example 1 to detect Hopf bifurcations. The zero contours of $R(s,\omega)$ and $I(s,\omega)$ with a viscosity contrast of 50 are shown in Figure \ref{fig:hopfContours_strat}. As before, each intersection of these curves indicates a Hopf bifurcation of frequency $\omega$ occurs at the $(Q_1,Q_C^*)$ pair associated with index $s$.  This figure is not symmetric because the underlying network is not. All the intersections are on the right side of the figure indicating that oscillations will occur when $Q_1$ is greater than $Q_2$. Unlike Figure \ref{fig:hopfContours}, here there is no differentiation between low frequency Hopf bifurcations that emerge near the saddle-node bifurcation and higher Hopf bifurcations that emerge away from it. All Hopf bifurcations in Figure \ref{fig:hopfContours_strat} emerge near the saddle-node bifurcation and grow towards the $Q_1 = 1$ boundary. At a viscosity contrast of 50 the $Q_1 = 1$ boundary has been destabilized. This fact is experimentally relevant, as oscillations would be observed with inlet 2 in Figure \ref{fig:schematic} shut off, leading to a simplified experimental design. As the viscosity contrast is increased, additional bands of instability appear and grow from the saddle-node bifurcation towards the $Q_1 = 1$ boundary.

\begin{figure}[!h]
\begin{center}
\includegraphics[width=90mm]{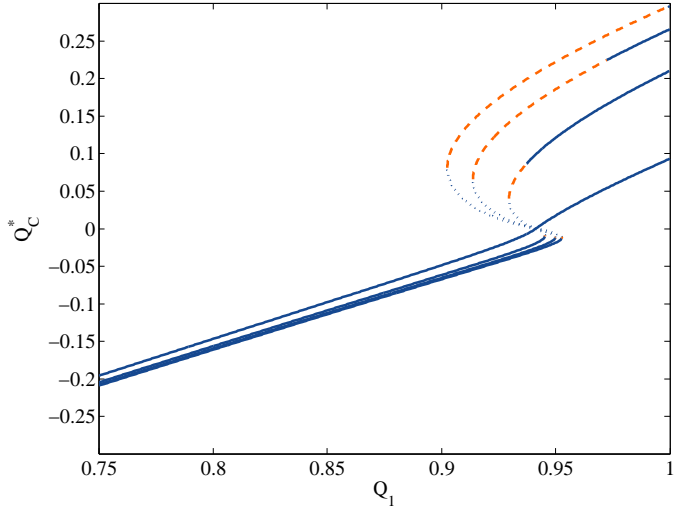}
\end{center}
\caption{Equilibrium curves of the stratified laminar flow model for viscosity contrast of 2, 10, 20, and 30. Solid regions of the equilibrium curve represent stable equilibria, while dotted regions represent unstable equilibria. Dashed regions indicate the existence of a limit cycle. Hopf bifurcations emerge from the points at which the equilibrium curves fold over. As the viscosity contrast is increased, the region of instability with $Q_C^* > 0$ grow towards the $Q_1 = 1$ boundary. }\label{fig:hopfHighlight_strat}
\end{figure}

In Figure \ref{fig:hopfHighlight_strat} we vary the viscosity contrast for a fixed geometry and track the associated bands of instability along the equilibrium curves. In this figure we plot the equilibrium curve in the $Q_1 \times Q_C^*$ plane at four values of the viscosity contrast. These curves are experimentally relevant as one can build a fixed network, change the inlet fluids to adjust the viscosity contrast and adjust the relative flow of the two inlets to move left and right along the $x$-axis \cite{karst2013}. Note that all the curves pass through the point $Q_C^*=0$ when $Q_1=16/17$. This trivial point is determined by noting that when $Q_C^*=0$ all the flow from inlet 1 goes through branch A and all the flow from inlet 2 goes through branch B. Since there is no flow in C, the pressure drop across A and B must be the same. Thus, the trivial point is given by
 \[
 Q_1 = \frac{  r_B \mu_2}{  r_B \mu_2 +  r_A \mu_1}.
 \]
 For our example $\mu_1=\mu_2$ and $\ell_A=\ell_B$, thus $Q_1=d_A^4/(d_B^4 + d_A^4)=16/17$.

The four curves shown in Figure \ref{fig:hopfHighlight_strat} represent the equilibrium solution for viscosity contrasts of 2, 10, 20, and 30. When the viscosity contrast is 2, the equilibrium curves are single-valued and there are no Hopf bifurcations. At a viscosity contrast of 10, the equilibrium curve becomes multi-valued over a small range around $Q_1=16/17$. For  this range of $Q_1$  there are two possible states, one with positive and negative $Q_C^*$. We also see a region of instability emerges  at the locations where the curves fold over. The instability band with positive $Q_C^*$ is much wider than the one with negative $Q_C^*$. As we increase the viscosity contrast to 20 the region of multiple equilibrium grows as does the band of instability. This behavior is different than the microvascular example in that the band of instability grows out of the point where the equilibrium curves fold over. We see that only the band with positive $Q_C^*$ grows significantly in size. When we increase the viscosity contrast to 30 the instability band encompasses the whole branch of the positive $Q_C^*$ equilibrium curve.

\begin{figure}[!h]
\centering
\includegraphics[width=130mm]{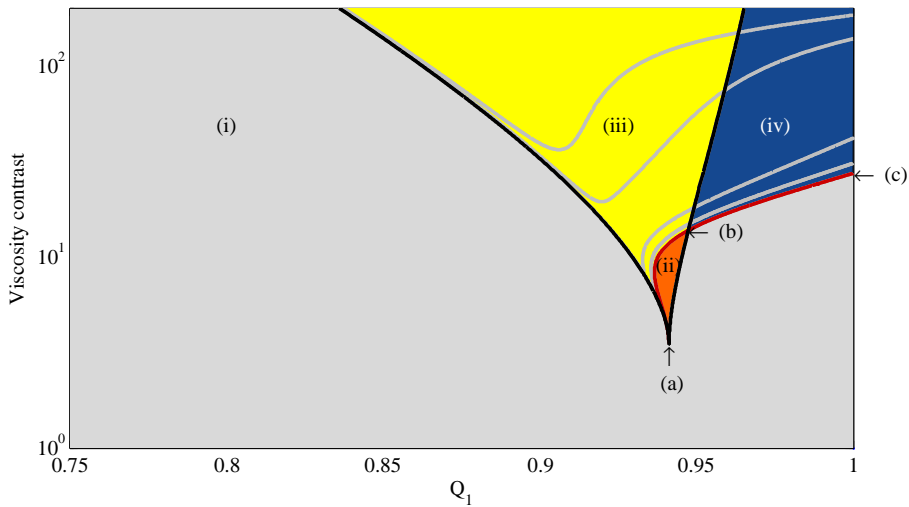}
\caption{Phase diagram in $Q_1 \times (\mu_\beta / \mu_\alpha)$ parameter space for the stratified laminar flow model. In the gray region (i), the system exhibits a single unique equilibrium state. In the orange region (ii), two stable equilibrium states exist. The yellow region (iii) represents parameters which support an unstable oscillation and one stable equilibrium state. The dark blue region (iv) represents parameters with a single oscillatory state.  The regions are separated by curves marking saddle-node bifurcations (black curves) and Hopf bifurcations (red/gray curves).}
\label{fig:bifurc_strat}
\end{figure}

We construct the phase diagram shown in Figure \ref{fig:bifurc_strat} to demonstrate the different possible states of the system. If we start with a low  viscosity contrast, \emph{i.e.}, less than $3.5$,  we have a single  equilibrium state for any $Q_1$. As we increase the viscosity contrast, multiple equilibrium behavior emerges from $Q_1=16/17$ when the viscosity contrast is $3.49$ (point a). As soon as multiple equilibrium exists, a Hopf bifurcation (denoted by the red curve) emerges from   the multiple equilibrium point. This instability occurs on the branch where $Q_C^*>0$; recall the behavior for a viscosity contrast of 10 from Figure \ref{fig:hopfHighlight_strat}. As we increase the viscosity contrast to $13.9$, this region  instability grows and eventually leaves the multiple equilibria region (point b).  After this viscosity contrast is exceeded we may have branches of the  equilibrium curve that are unstable via Hopf bifurcation, and there is no other possible stable equilibrium state. As we increase the viscosity contrast to $27.6$, the the region of instability reaches the boundary where $Q_1=1$ (point c); thus instability encompasses the entire branch of the equilibrium curve where $Q_C^*>0$.

A very narrow band of instability also grows along the right edge of the multiple equilibrium boundary. This band corresponds to the instability region seen for negative $Q_C^*$ at the fold in the equilibrium curve in Figure \ref{fig:hopfHighlight_strat}. This region is so narrow and only exists right the multiple equilibrium boundary that is likely of little practical interest and not observable. Due to the broken symmetry for this parameter set, we only see significant instability for cases where $Q_C^*>0$. Thus unlike the example with microvascular blood flow, here we do not find co-existing limit cycles.

\begin{figure}[!h]
\begin{center}
\includegraphics[width=125mm]{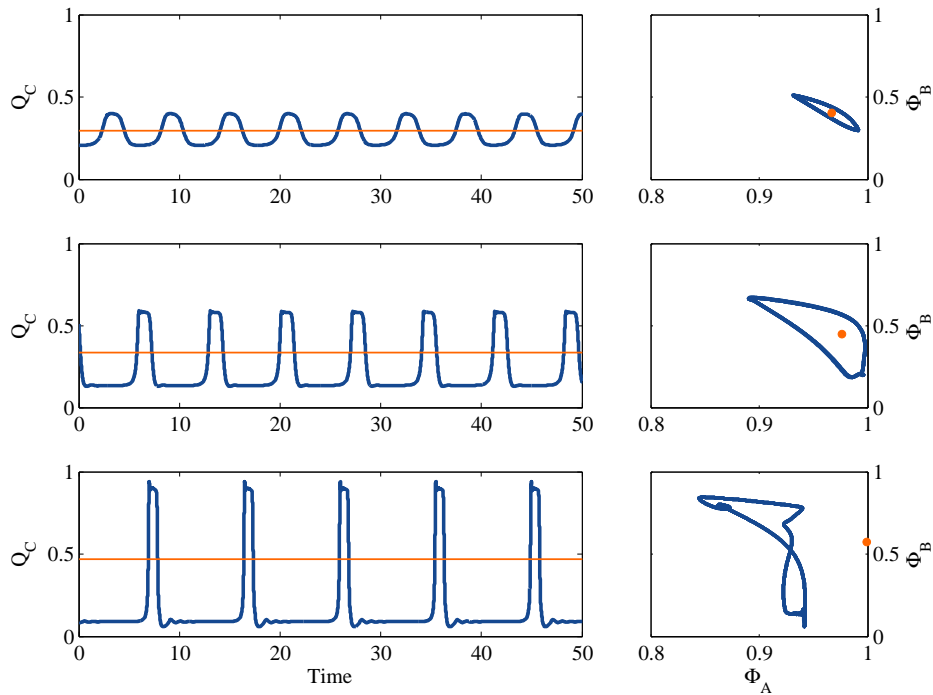}
\end{center}
\caption{Time series and associated phase plot from the direct numerical simulation of the stratified flow model with $Q_1 = 1$ for different values of viscosity contrast, 30, 50, and 500 from top to bottom. In each phase plot the equilibrium solution is shown as the dot. As the viscosity contrast in increased, progressively higher frequency bands of instability reach the $Q_1 = 1$ boundary in Figure \ref{fig:bifurc_strat}. The presence of these higher frequencies result in richer temporal dynamics in $Q_C(t)$ at higher viscosity contrasts.}\label{fig:timeSeries_strat}
\end{figure}

As in Example 1, while the linear analysis can provide some insight into the types of behaviors we may see, we must resort to full numerical simulation in order to see the complete dynamics. In Figure \ref{fig:timeSeries_strat} we show some sample dynamics for the stratified flow model with $Q_1 = 1$. With the viscosity contrast set to 30, we observe a relatively sinusoidal oscillation in $Q_C(t)$. As the viscosity contrast is increased, Figure \ref{fig:bifurc_strat} shows that higher frequency bands of instability grow towards the $Q_1 = 1$ boundary. With the viscosity contrast set to 50, for instance, there are 3 distinct bands of instability that have crossed the $Q_1 = 1$ boundary. These additional frequencies lead to richer temporal dynamics in the flow $Q_C(t)$ as seen in the middle pane of Figure \ref{fig:timeSeries_strat}. As the viscosity contrast is increased to 500, progressively higher frequency Hopf bifurcations have crossed the $Q_1 = 1$ boundary. This broader spectrum manifests as abrupt changes in the flow rate $Q_C(t)$ as seen in the bottom pane of Figure \ref{fig:timeSeries_strat}.

\section{Conclusions}
We have demonstrated a rich set of dynamics which emerge from  simple fluid networks with practical and experimental relevance. We have presented a method for analyzing these fluid networks which has a large number of important free parameters. Through direct numerical simulation, the parameter space is too large to span in a systematic way. We find large ranges of parameter space in which equilibrium solutions to the phase and flow distribution within a network are unstable and  spontaneous oscillations may emerge. We also find complex nonlinear dynamics for large viscosity contrasts.

While we have presented our results in a manner which is experimentally relevant, the details of the constitutive laws are such that they are critical to the exact predictions of stability and are difficult to experimentally control. Thus while our laws for viscosity and phase separation at a node are realistic for blood flow, the viscosity contrast of blood (contrast between plasma and red cell rich fluid) is limited to approximately 10, thus the contrast of 30 or 50 to see oscillations is probably still out of experimental range. However, through careful selection of the network parameters it may be possible to find examples which occur in realistic experimental systems. Further, the range of parameters where spontaneous oscillations exist for this network is much broader and more realistic than the equivalent 2-node network~\cite{Geddes:2007}, thus adding an additional network branch might be sufficient to bring the
dynamics into experimental space.

On the other hand, the predictions for the stratified network model are well within the range of what is possible experimentally~\cite{karst2013}. The stratified system has the advantage that viscosity is a more easily controlled parameter through the selection of the fluids and the flow state is a natural consequence of buoyancy effects.  On going work is aimed at direct observation of these predictions.

\section*{Acknowledgments}
This work was supported in part by the National Science Foundation under Contract No. DMS-1211640.

\bibliographystyle{siam}

\end{document}